\documentclass[12pt]{article}      % Specifies the document class \usepackage{epsfig}

\title{Bernoulli-Fibonacci Polynomials}

\author{ Oktay K Pashaev and Merve $\rm{\ddot{O}}$zvatan\\Department of Mathematics\\ Izmir Institute of Technology \\ Urla-Izmir, 35430, Turkey}

\begin{document}
\newcommand{\be}{\begin{equation}}
\newcommand{\ee}{\end{equation}}
\newcommand{\bea}{\begin{eqnarray}}
\newcommand{\eea}{\end{eqnarray}}
\newcommand{\disp}{\displaystyle}
\newcommand{\la}{\langle}
\newcommand{\ra}{\rangle}

\newtheorem{thm}{Theorem}[subsection]
\newtheorem{cor}[thm]{Corollary}
\newtheorem{lem}[thm]{Lemma}
\newtheorem{prop}[thm]{Proposition}
\newtheorem{definition}[thm]{Definition}
\newtheorem{rem}[thm]{Remark}
\newtheorem{prf}[thm]{Proof}

\maketitle

%%%%%%%%%%%%%%%%%%%%%%%%%%%%%%%%%%%%%%%%%%%%%%%%%%%%%%%%%%%%%% % You may repeat \author \address as often as necessary      % %%%%%%%%%%%%%%%%%%%%%%%%%%%%%%%%%%%%%%%%%%%%%%%%%%%%%%%%%%%%%%

\begin{abstract}
By using definition of Golden derivative, corresponding Golden exponential function and Fibonomial coefficients, we introduce generating
functions for Bernoulli-Fibonacci polynomials and related numbers. Properties of these polynomials and numbers are studied in parallel with
usual Bernoulli counterparts. Fibonacci numbers and Golden ratio are intrinsically involved in formulas obtained.
\end{abstract}

\section{Introduction}

The Bernoulli numbers and related polynomials can be introduced on the basis of the so called h-calculus \cite{Kac}. Several
generalizations of these numbers in q-calculus are known and listed in papers \cite{qBern}, \cite{mourad}. Here we are going to introduce
analogs of Bernoulli polynomials in the so called Golden calculus, which is based on the Golden ratio dilatations and Binet formula
for Fibonacci numbers \cite{golden}.

The paper is organized in following way. Firstly, in Section 2 we briefly review main properties of Bernoulli polynomials
. In Section 3 we introduce main objects of Golden calculus, as Golden derivative, exponential function and Fibonomials.
In the main part, Section 4, we introduce Bernoulli-Fibonacci polynomials and corresponding numbers, formulate and proof related theorems.

\section{Bernoulli polynomials}

Here we review briefly main definitions and properties of Bernoulli polynomials and Bernoulli numbers \cite{Kac}.

\begin{definition}
The generating function for Bernoulli polynomials is defined by Taylor series expansion,
\begin{equation}
 \frac{z \phantom{a}e^{zx}}{e^{z}-1}=\sum_{n=0}^\infty B_{n}(x) \frac{z^n}{n!}, \label{genfunctionforb.p}
\end{equation}
where $B_{n}(x)$ are the Bernoulli polynomials in x, for all $n>0$.

Bernoulli numbers are a special values of the Bernoulli polynomials $B_n(x)$, $b_n=B_n(0)$. The generating function for Bernoulli numbers is
\bea
\frac{z}{e^{z}-1}=\sum_{n=0}^\infty b_{n} \frac{z^n}{{n}!}. \label{genfunctionforb.n}
\eea
\end{definition}
First few Bernoulli polynomials 
\begin{eqnarray}
B_{0}(x)&=&1 \nonumber \\
B_{1}(x)&=&x-\frac{1}{2}\nonumber \\
B_{2}(x)&=&x^2-x+\frac{1}{6}	 \nonumber \\
B_{3}(x)&=&x^3-\frac{3}{2} x^2+\frac{1}{2} x	 \nonumber \\
B_{4}(x)&=&x^4-2x^3+x^2-\frac{1}{30}	\nonumber \\
B_{5}(x)&=&x^5-\frac{5}{2} x^4+\frac{5}{3} x^3-\frac{1}{6}x	\nonumber \\
B_{6}(x)&=&x^6-3x^5+\frac{5}{2} x^4-\frac{1}{2}x^2+\frac{1}{42}, \nonumber
\end{eqnarray}
determine corresponding Bernoulli numbers \\
$$b_{0}=1, \quad b_{1}=-\frac{1}{2},\quad b_{2}=\frac{1}{6},\quad b_{3}=0,\quad b_{4}=-\frac{1}{30},\quad b_{5}=0,\quad
b_{6}=\frac{1}{42}.$$
 Bernoulli polynomials and numbers poses following properties\cite{Kac}.

\begin{prop}
For odd Bernoulli numbers
\bea
b_{2n+1}=0,
\eea
where $n=1,2,\ldots$
\end{prop}

\begin{prop} \label{prop3.2}Derivative of $n^{th}$ Bernoulli polynomial gives,
\bea
\frac{d}{dx}B_{n}(x)=B^{'}_{n}(x)=n\phantom{.} B_{n}(x).
\eea
\end{prop}

\begin{prop}\label{prop3.3}For any $n \geq 0$, we have expansion
\bea
B_{n}(x)=\sum_{j=0}^n {n \choose j}\phantom{.}  b_{j}\phantom{a} x^{n-j}.
\eea
\end{prop}

\begin{prop}\label{prop3.4}For any $n \geq 1$, we have summation formula
\bea
\sum_{j=0}^{n-1} {n \choose j} \phantom{.} B_{j} (x) =n \phantom{a} x^{n-1}.  \nonumber
\eea
\end{prop}

\begin{cor}\label{corollary3.1} From previous proposition at $x=0$, for any $n \geq 2$ we have
\bea
\sum_{j=0}^{n-1}\phantom{.} {n \choose j} \phantom{.} b_{j}=0.  \nonumber
\eea
\end{cor}

\begin{prop}\label{prop3.5}For any $n \geq 2$,
\bea
B_{n}(1)=b_{n}.
\eea
\end{prop}

\section{Golden Derivative and Golden Exponential}

The main ingredients of Golden calculus are described in \cite{golden}. It is based on Fibonacci numbers
and Binet representation of these numbers:
\be  F_n = \frac{\varphi ^n- {\varphi'}^n}{\varphi- \varphi'},\label{fbbinet}\ee where $\varphi, \, \varphi'$ are positive and negative roots of equation $x^2-x-1=0.$
These roots  explicitly are \be  \varphi= \frac{1+\sqrt{5}}{2}, \,\,\,\,\,\,\,\,\, \varphi'= \frac{1-\sqrt{5}}{2}= -\frac{1}{\varphi}. \ee
Number $\varphi$ is known as the Golden ratio or the Golden section.
The Golden derivative operator is defined as
\be  D_F f(x)= \frac{f(\varphi x)-f(-\frac{x}{\varphi})}{(\varphi+\frac{1}{\varphi})x} \label{goldenderivative} .\ee
For analytic function $f(x)$ it can be represented by Binet formula for dilatation  operator

\be  F_{x \frac{d}{dx}} = \frac{\varphi^{x \frac{d}{dx}}- \varphi'^{x \frac{d}{dx}}}{\varphi - \varphi'} , \label{gderivative}\ee
so that $$x F_{x \frac{d}{dx}}f(x) = D_F f(x).$$
The Golden exponential function is entire analytic function, defined as
\be e_F^x\equiv \sum_{n=0}^\infty \frac{x^n}{F_n!},\label{exp}\ee
 and satisfying
\be D_F e_F^{k x}= k e_F^{k x}.\label{DeF}\ee
The  Golden binomial, introduced in \cite{golden}, can be expanded
\bea (x+y)_F^n \equiv (x+y)_{\varphi,-\frac{1}{\varphi}}^n &=&\sum^{n}_{k=0}{ n \brack k}_{F} (-1)^{\frac{k(k-1)}{2}} x^{n-k} y^k  \label{goldenbinomexpansion}\eea
in terms of the Golden binomial coefficients, called as Fibonomials,
\be { n \brack k}_F=  \frac{F_n!}{F_{n-k}! F_k!}, \label{goldenbinom}\ee
with $n$ and $k$ being nonnegative integers, $n\geq k$ .
  They satisfy next recursion formulas
\bea { n \brack k}_{F}&=& \left(-\frac{1}{\varphi}\right)^k { n-1 \brack k}_{F} + \varphi^{n-k} { n-1 \brack k-1}_{F} \label{goldenpascal1}\\
&=& \varphi^k { n-1 \brack k}_{F} + \left(-\frac{1}{\varphi}\right)^{n-k} { n-1 \brack k-1}_{F} .\label{goldenpascal2} \eea

\section{Bernoulli-Fibonacci polynomials}

 By using definition of Golden exponential function $e^{x}_{F}$ in $(\ref{exp})$ we introduce generating function
 for Bernoulli-Fibonacci polynomials.

\begin{definition}
Generating function for Bernoulli-Fibonacci polynomials $B_{n}^{F}(x)$ is defined by series expansion,
\begin{eqnarray}
 \frac{z \phantom{a}e_{F}^{zx}}{e_{F}^{z}-1}=\sum_{n=0}^\infty B_{n}^{F}(x) \frac{z^n}{F_{n}!}.  \label{generatingfunctionforbernoullif.p}
\end{eqnarray}
The Bernoulli-Fibonacci numbers are special values of these polynomials,
\bea
b_{n}^{F}=B_{n}^{F}(0),
\eea
with generating function,
\bea
\frac{z}{e_{F}^{z}-1}=\sum_{n=0}^\infty b_{n}^{F} \frac{z^n}{F_{n}!}. \label{generatingfunctionforbernoullif.n}
\eea
\end{definition}
For the first few polynomials we have,
\begin{eqnarray}
B_{0}^{F}(x)&=&1 ,\nonumber \\
B_{1}^{F}(x)&=&\frac{x}{F_{1}!}-\frac{1}{F_{2}!} ,\nonumber \\
B_{2}^{F}(x)&=&x^2-x\frac{1}{F_{1}!}+\frac{1}{F_{2}!}-\frac{1}{F_{3}}, \nonumber \\
B_{3}^{F}(x)&=&x^3-x^2\frac{F_{3}!}{(F_{2}!)^2}+x\left(\frac{F_{3}!}{F_{1}!(F_{2}!)^2}-\frac{1}{F_{1}!}\right)+\frac{2}{F_{2}!}-\frac{1}{F_{4}}-F_{3}, \nonumber \\
B_{4}^{F}(x)&=&x^4- x^3 \left(\frac{F_{4}!}{F_{3}!F_{2}!}\right)+ x^2 \left(\frac{F_{4}!}{(F_{2}!)^3}-\frac{F_{4}!}{F_{2}! F_{3}!}\right) \nonumber \\
&+&x \left(-\frac{F_{4}!}{F_{1}! (F_{2}!)^3}+2 \frac{F_{4}!}{F_{1}! F_{2}! F_{3}!}-\frac{1}{F_{1}!} \right) \nonumber \\
&+& \left(\frac{F_{4}!}{ (F_{2}!)^4}-3 \frac{F_{4}!}{(F_{2}!)^2 F_{3}!}+2\frac{1}{F_{2}!}+\frac{F_{4}!}{(F_{3}!)^2}-\frac{1}{F_{5}} \right),\nonumber \\
B_{5}^{F}(x)&=&x^5+x^4 \left( -\frac{F_{5}!}{F_{4}! F_{2}!} \right)+ x^3 \left( -\frac{F_{5}!}{(F_{3}!)^2}+ \frac{F_{5}!}{F_{3}! (F_{2}!)^2}\right) \nonumber \\
&+&   x^2 \left( -\frac{F_{5}!}{F_{2}! F_{4}!}+2 \frac{F_{5}!}{(F_{2}!)^2 F_{3}!}-\frac{F_{5}!}{(F_{2}!)^4} \right)  \nonumber \\
&+&   x \left( -\frac{1}{F_{1}!}+ \frac{F_{5}!}{F_{1}!(F_{3}!)^2}+2\frac{F_{5}!}{F_{1}!F_{2}!F_{4}!}-3\frac{F_{5}!}{F_{1}!(F_{2}!)^2F_{3}!}+\frac{F_{5}!}{F_{1}!(F_{2}!)^4} \right)  \nonumber \\
&+&  \left(-\frac{F_{5}!}{(F_{2}!)^5}+4 \frac{F_{5}!}{(F_{2}!)^3 F_{3}!}-3 \frac{F_{5}!}{(F_{3}!)^2 F_{2}!}-3 \frac{F_{5}!}{(F_{2}!)^2 F_{4}!}+2 \frac{F_{5}!}{F_{3}! F_{4}!}+\frac{2}{F_{2}!}- \frac{F_{5}!}{F_{6}!}   \right) \nonumber
\end{eqnarray}

It is noticed that all coefficients of Bernoulli-Fibonacci polynomials are rational Fibonacci numbers.

\begin{prop}(Compare with Proposition $(\ref{prop3.2})$)

The Golden derivative, applied to Bernoulli-Fibonacci polynomials $B_{n}^{F}(x)$ gives Fibonacci numbers,
\bea
D^{x}_{F}(B_{n}^{F}(x))= F_{n} \phantom{.} B_{n-1}^{F}(x). \label{goldenderivativeapplicationtobernoullipolynomials}
\eea
\end{prop}
\begin{prf}
Taking the derivative from both sides of equation $(\ref{generatingfunctionforbernoullif.p})$ gives,
\bea
D^{x}_{F}\left( \frac{z \phantom{a}e_{F}^{zx}}{e_{F}^{z}-1}\right)&=&D^{x}_{F}\left(  \sum_{n=0}^\infty B_{n}^{F}(x) \frac{z^n}{F_{n}!} \right) \nonumber
\Rightarrow
\\
\frac{z \phantom{a}D^{x}_{F}(e_{F}^{zx})}{e_{F}^{z}-1}&=&D^{x}_{F}\left( B_{0}^{F}(x)+B_{1}^{F}(x)\frac{z}{F_{1}!}+B_{2}^{F}(x)\frac{z^2}{F_{2}!}+ \ldots \right) \nonumber
\eea
In the l.h.s., Golden derivative can be calculated from equation $(\ref{DeF})$. In the r.h.s. we use
\bea
D^{x}_{F} (B_{0}^{F}(x))&=& D^{x}_{F} (1)=0 .\nonumber
\eea
Then,
\bea
z \cdot \frac{z e_{F}^{zx}}{e_{F}^{z}-1}&=&  \sum_{k=1}^\infty D^{x}_{F}\left(B_{k}^{F}(x)\right) \frac{z^k}{F_{k}!} \Rightarrow \nonumber \\
z \sum_{n=0}^\infty B_{n}^{F}(x) \frac{z^n}{F_{n}!} &=&  \sum_{k=1}^\infty D^{x}_{F}\left(B_{k}^{F}(x)\right) \frac{z^k}{F_{k}!} \Rightarrow \nonumber \\
\sum_{n=0}^\infty B_{n}^{F}(x) \frac{z^{n+1}}{F_{n}!} &=&  \sum_{k=1}^\infty D^{x}_{F}\left(B_{k}^{F}(x)\right) \frac{z^k}{F_{k}!}. \nonumber
\eea
In the r.h.s., by shifting $k-1=n$ we have
\bea
\sum_{n=0}^\infty B_{n}^{F}(x) \frac{z^{n+1}}{F_{n}!} &=&  \sum_{n=0}^\infty D^{x}_{F}\left(B_{n+1}^{F}(x)\right) \frac{z^{n+1}}{F_{{n+1}}!} \Rightarrow \nonumber \\
\sum_{n=0}^\infty B_{n}^{F}(x) \frac{z^{n+1}}{F_{n}!} &=&  \sum_{n=0}^\infty D^{x}_{F}\left(B_{n+1}^{F}(x)\right) \frac{1}{F_{n+1}}\frac{z^{n+1}}{F_{n}!}. \nonumber
\eea
From equality of these two series, finally we  get
\bea
B_{n}^{F}(x)= \frac{D^{x}_{F}\left(B_{n+1}^{F}(x)\right)}{F_{n+1}}. \nonumber
\eea
\end{prf}

From Bernoulli-Fibonacci polynomials, due to $b_{n}^{F}=B_{n}^{F}(0)$, we find Bernoulli-Fibonacci numbers  as,
\bea
b_{0}^{F}=1, \phantom{..} b_{1}^{F}=-1,\phantom{..} b_{2}^{F}=\frac{1}{2},\phantom{..} b_{3}^{F}=-\frac{1}{3},\phantom{..} b_{4}^{F}=\frac{3}{10},\phantom{..} b_{5}^{F}=-\frac{5}{8},\phantom{..}b_{6}^{F}=\frac{101}{39},... \nonumber
\eea
\begin{prop}(Compare with Proposition $(\ref{prop3.3})$)

Bernoulli-Fibonacci Polynomials can be represented by sum,
\begin{eqnarray}
B_{n}^{F} (x)= \sum_{j=0}^n {n \brack j}_{F} b_{j}^{F}\phantom{a} x^{n-j}. \label{alternativerepresentationofb.f.p.}
\end{eqnarray}
\end{prop}
\begin{prf}
It is evident from this formula that $B_{n}^{F} (0) = b_{n}^{F}$. So it suffices to show that
 $B_{n}^{F}(x)$  are satisfying equation (\ref{goldenderivativeapplicationtobernoullipolynomials}):
\bea
D^{x}_{F} \left[ B_{n}^{F} (x) \right]&{(\ref{alternativerepresentationofb.f.p.})}{=}&D^{x}_{F} \left[ \sum_{j=0}^n {n \brack j}_{F} b_{j}^{F}\phantom{a} x^{n-j} \right]  \nonumber \\
&=& \sum_{j=0}^{n-1} {n \brack j}_{F} b_{j}^{F}\phantom{a} D^{x}_{F}\left( x^{n-j} \right) \nonumber \\
&=& \sum_{j=0}^{n-1} {n \brack j}_{F} b_{j}^{F}\phantom{a}F_{n-j}\phantom{a}  x^{n-j-1} \nonumber \\
&=& \sum_{j=0}^{n-1} \frac{F_{n}!}{F_{n-j}!\phantom{a}F_{j}!}\, b_{j}^{F}\phantom{a}F_{n-j}\phantom{a}  x^{n-j-1} \nonumber \\
&=& \sum_{j=0}^{n-1} \frac{F_{n}!}{F_{n-j-1}!\phantom{a}F_{j}!}\, b_{j}^{F}\phantom{a}  x^{n-j-1} \nonumber \\
&=& F_{n} \sum_{j=0}^{n-1} \frac{F_{n-1}!}{F_{n-j-1}!\phantom{a}F_{j}!}\, b_{j}^{F}\phantom{a}  x^{n-j-1} \nonumber \\
&{(\ref{alternativerepresentationofb.f.p.})}{=}& F_{n} \phantom{a} B_{n-1}^{F}(x). \nonumber
\eea
\end{prf}

\begin{prop}(Compare with Proposition $(\ref{prop3.4})$)

For $n \geq 1$, Bernoulli-Fibonacci polynomials can be calculated recursively by formula:
\bea
\sum_{l=0}^{n-1} {n \brack l}_{F} B_{l}^{F} (x) =F_{n} \phantom{a} x^{n-1}.  \nonumber
\eea
\end{prop}
\begin{prf} Starting from the generating function,
\bea
\frac{z \phantom{a}e_{F}^{zx}}{e_{F}^{z}-1}=\sum_{n=0}^\infty B_{n}^{F}(x) \frac{z^n}{F_{n}!}, \nonumber
\eea
multiplying it by $e_{F}^{z}$;
\bea
\frac{z \phantom{a}e_{F}^{zx}\phantom{a}e_{F}^{z}}{e_{F}^{z}-1}=\sum_{n=0}^\infty B_{n}^{F}(x)\phantom{a} e_{F}^{z} \frac{z^n}{F_{n}!},  \nonumber
\eea
and taking difference, we get
\bea
\frac{z \phantom{a}e_{F}^{zx}}{e_{F}^{z}-1}\left(e_{F}^{z}-1 \right)=\sum_{n=0}^\infty\left(B_{n}^{F}(x) \phantom{a}e_{F}^{z}- B_{n}^{F}(x) \right)\frac{z^n}{F_{n}!} \Rightarrow \nonumber \\
z \phantom{a}e_{F}^{zx}=\sum_{n=0}^\infty\left( B_{n}^{F}(x) \phantom{a}e_{F}^{z}- B_{n}^{F}(x) \right)\frac{z^n}{F_{n}!} \Rightarrow \nonumber \\
D^{x}_{F}\left( e_{F}^{zx}\right)=\sum_{n=0}^\infty\left( B_{n}^{F}(x) \phantom{a}e_{F}^{z}- B_{n}^{F}(x) \right)\frac{z^n}{F_{n}!} \Rightarrow  \nonumber \\
D^{x}_{F}\left(\sum_{n=0}^\infty \frac{\left( z x \right)^n}{F_{n}!} \right)=\sum_{n=0}^\infty \left( B_{n}^{F}(x) \phantom{a}e_{F}^{z}- B_{n}^{F}(x) \right)\frac{z^n}{F_{n}!} \Rightarrow \nonumber \\
D^{x}_{F}\left(\sum_{n=1}^\infty \frac{ z^n \phantom{.} x^n}{F_{n}!} \right)=\sum_{n=0}^\infty \left( B_{n}^{F}(x) \phantom{a}e_{F}^{z}- B_{n}^{F}(x) \right)\frac{z^n}{F_{n}!} \Rightarrow \nonumber \\
\sum_{n=1}^\infty \frac{ z^n \phantom{.}D^{x}_{F} \phantom{.} \left( x^{n}\right)}{F_{n}!} =\sum_{n=0}^\infty \left( B_{n}^{F}(x) \phantom{a}e_{F}^{z}- B_{n}^{F}(x) \right)\frac{z^n}{F_{n}!} \Rightarrow \nonumber \\
\sum_{n=1}^\infty \frac{ z^n \phantom{.}F_{n}\phantom{.} x^{n-1}}{F_{n}!}=\sum_{l=0}^\infty  B_{l}^{F}(x) \phantom{a}e_{F}^{z} \frac{z^l}{F_{l}!}-\sum_{n=0}^\infty  B_{n}^{F}(x) \frac{z^n}{F_{n}!}. \nonumber
\eea
For the first sum in r.h.s. of this equation,
\bea
\sum_{l=0}^\infty  B_{l}^{F}(x) \phantom{a}e_{F}^{z} \frac{z^l}{F_{l}!}&=&\sum_{l=0}^\infty \sum_{k=0}^\infty  B_{l}^{F}(x) \frac{z^k}{F_{k}!} \frac{z^l}{F_{l}!}=\sum_{l=0}^\infty \sum_{k=0}^\infty  B_{l}^{F}(x) \frac{z^{k+l}}{F_{k}! F_{l}!} \nonumber
\eea
\bea
&=&\sum_{l=0}^\infty \sum_{k=0}^\infty  B_{l}^{F}(x) \frac{z^{k+l}}{F_{k}! F_{l}!} \nonumber \\
&{k+l=n}{=}&\sum_{n=0}^\infty \frac{1}{F_{n}!} \sum_{l=0}^n  B_{l}^{F}(x) \frac{z^{n}\phantom{.}F_{n}!  }{F_{n-l}! F_{l}!} \nonumber \\
&=&\sum_{n=0}^\infty \frac{z^{n}}{F_{n}!} \left(  \sum_{l=0}^n  {n \brack l}_{F} B_{l}^{F}(x) \right). \nonumber
\eea
After substituting this,
\bea
\sum_{n=1}^\infty \frac{ z^n \phantom{.}F_{n}\phantom{.} x^{n-1}}{F_{n}!}&=&\sum_{n=0}^\infty \frac{z^{n}}{F_{n}!} \left(  \sum_{l=0}^n  {n \brack l}_{F} B_{l}^{F}(x) \right)-\sum_{n=0}^\infty  B_{n}^{F}(x) \frac{z^n}{F_{n}!} \Rightarrow \nonumber \\
\sum_{n=1}^\infty \frac{ z^n \phantom{.}F_{n}\phantom{.} x^{n-1}}{F_{n}!}&=& B_{0}^{F}(x) {0 \brack 0}_{F}+ \sum_{n=1}^\infty \frac{z^{n}}{F_{n}!} \left(  \sum_{l=0}^n  {n \brack l}_{F} B_{l}^{F}(x) \right)-B_{0}^{F}(x)-\sum_{n=1}^\infty  B_{n}^{F}(x) \frac{z^n}{F_{n}!} \Rightarrow  \nonumber \\
\sum_{n=1}^\infty \frac{ z^n}{F_{n}!}  \phantom{.}F_{n}\phantom{.} x^{n-1}&=& \sum_{n=1}^\infty \frac{z^{n}}{F_{n}!} \left(  \sum_{l=0}^n  {n \brack l}_{F} B_{l}^{F}(x) - B_{n}^{F}(x)\right). \nonumber
\eea
Then, we have
\bea
\sum_{l=0}^n  {n \brack l}_{F} B_{l}^{F}(x)-B_{n}^{F}(x)=F_{n}\phantom{.} x^{n-1}, \nonumber
\eea
where $n\geq1$. Cancelling  $n^{th}$ terms,
\bea
\sum_{l=0}^{n-1}  {n \brack l}_{F} B_{l}^{F}(x)+B_{n}^{F}(x)\phantom{a} {n \brack n}_{F} -B_{n}^{F}(x)=F_{n}\phantom{.} x^{n-1}, \nonumber
\eea
finally, we get desired equality,
\bea
\sum_{l=0}^{n-1}  {n \brack l}_{F} B_{l}^{F}(x)=F_{n}\phantom{.} x^{n-1}. \nonumber
\eea
\end{prf}
\begin{cor}(Compare with Corollary $(\ref{corollary3.1})$)

From previous proposition, by taking $x=0$, we have the formula, allowing us to compute Bernoulli-Fibonacci numbers inductively,
\bea
\sum_{j=0}^{n-1} {n \brack j}_{F} b_{j}^{F}=0,
\eea
where $n \geq 2$.
\end{cor}
\begin{prf}
Proof will be done by equating expansions of $B_{n}^{F}(x)+F_n x^{n-1}$, obtained in two ways. By the first way we obtain,
\begin{eqnarray}
B_{n}^{F}(x)+F_n x^{n-1}=x^n + \sum_{j=2}^{n} {n \brack j}_{F} b_{j}^{F} x^{n-j}, \label{ononehand}
\end{eqnarray}
and by the second one,
\begin{eqnarray}
B_{n}^{F}(x)+F_n x^{n-1} \equiv H_{n}(x)=\sum_{k=0}^{n} {n \brack k}_{F} B_{n-k}^{F}(x). \label{ontheotherhand}
\end{eqnarray}
These equations are derived in Appendix.
Comparison of these two expansions gives,
\begin{eqnarray}
\sum_{k=0}^{n} {n \brack k}_{F} B_{n-k}^{F}(x)=x^n + \sum_{j=2}^{n} {n \brack j}_{F} b_{j}^{F} x^{n-j}. \nonumber
\end{eqnarray}
Then we have,
\bea
\sum_{k=0}^{n-1} {n \brack k}_{F} B_{n-k}^{F}(x)+B_{0}^{F}(x)=x^n + \sum_{j=2}^{n} {n \brack j}_{F} b_{j}^{F} x^{n-j} \Rightarrow \nonumber
\eea
\bea
x^n-1=\sum_{k=0}^{n-1} {n \brack k}_{F} B_{n-k}^{F}(x)- \sum_{j=2}^{n} {n \brack j}_{F} b_{j}^{F} x^{n-j} \Rightarrow  \nonumber
\eea
\bea
x^n-1=\sum_{k=0}^{n-1} {n \brack k}_{F} B_{n-k}^{F}(x)- \left( {n \brack 2}_{F} b_{2}^{F} x^{n-2}+\ldots+{n \brack n-1}_{F} b_{n-1}^{F} x+{n \brack n}_{F} b_{n}^{F} \right). \nonumber
\eea
For $x=0$ it gives,
\begin{eqnarray}
-1=\sum_{k=0}^{n-1} {n \brack k}_{F} b_{n-k}^{F} - b_{n}^{F}. \nonumber
\end{eqnarray}
After simplifying,
\bea
-1=\sum_{k=1}^{n-1} {n \brack k}_{F} b_{n-k}^{F}+  {n \brack 0}_{F} b_{n}^{F} - b_{n}^{F}, \nonumber
\eea
 we get
\bea
\sum_{k=1}^{n-1} {n \brack k}_{F} b_{n-k}^{F}+1=0. \nonumber
\eea
By using symmetry of Fibonomial coefficients, we  rewrite it as,
\bea
\sum_{k=1}^{n-1} {n \brack n-k}_{F} b_{n-k}^{F}+1=0,   \nonumber
\eea
then by denoting $n-k=j$,
\bea
\sum_{j=n-1}^{1} {n \brack j}_{F} b_{j}^{F}+1=0.   \nonumber
\eea
Substituting $1 = b_{0}^{F}$,
\bea
\sum_{j=1}^{n-1} {n \brack j}_{F} b_{j}^{F}+{n \brack 0}_{F} b_{0}^{F} =0 \nonumber
\eea
 we obtain desired result
\bea
\sum_{j=0}^{n-1} {n \brack j}_{F} b_{j}^{F}=0 .\nonumber
\eea
\end{prf}

\begin{prop}(Compare with Proposition $(\ref{prop3.5})$)

For any $n \geq 2$, Bernoulli-Fibonacci polynomials and corresponding numbers satisfy the following equation,
\bea
 B_{n}^{F}(1)=b_{n}^{F}.
\eea

\end{prop}
\begin{prf}Starting with,
\bea
\frac{e_{F}^{zx}-1}{e_{F}^{z}-1}&=&\left(\frac{e_{F}^{zx}}{e_{F}^{z}-1}-\frac{1}{e_{F}^{z}-1}\right)=\left(\frac{z \phantom{a}e_{F}^{zx}}{e_{F}^{z}-1}-\frac{z}{e_{F}^{z}-1}\right)\frac{1}{z} \nonumber \\
&{(\ref{generatingfunctionforbernoullif.p})}{=}& \left(\sum_{n=0}^\infty B_{n}^{F}(x) \frac{z^n}{F_{n}!}-\frac{z}{e_{F}^{z}-1}\right)\frac{1}{z} \nonumber \\
&{(\ref{generatingfunctionforbernoullif.n})}{=}&  \left(\sum_{n=0}^\infty B_{n}^{F}(x) \frac{z^n}{F_{n}!}-\sum_{n=0}^\infty b_{n}^{F} \frac{z^n}{F_{n}!}\right)\frac{1}{z} \nonumber \\
&=&\sum_{n=0}^\infty \left( B_{n}^{F}(x)-  b_{n}^{F}\right) \frac{z^{n-1}}{F_{n}!}. \nonumber
\eea

Thus, we have,
\bea
\frac{e_{F}^{zx}-1}{e_{F}^{z}-1}&=&\left(B_{0}^{F}(x)-b_{0}^{F} \right) \frac{z^{-1}}{F_{0}!}+\left(B_{1}^{F}(x)-b_{1}^{F} \right) \frac{1}{F_{1}!}+\sum_{n=2}^\infty \left( B_{n}^{F}(x)-  b_{n}^{F}\right) \frac{z^{n-1}}{F_{n}!} \Rightarrow  \nonumber \\
\frac{e_{F}^{zx}-1}{e_{F}^{z}-1}&=&\left(1-1 \right) \frac{z^{-1}}{F_{0}!}+\left(x-1-(-1) \right) \frac{1}{F_{1}!}+\sum_{n=2}^\infty \left( B_{n}^{F}(x)-  b_{n}^{F}\right) \frac{z^{n-1}}{F_{n}!} \Rightarrow \nonumber \\
\frac{e_{F}^{zx}-1}{e_{F}^{z}-1}&=&x+\sum_{n=2}^\infty \left( B_{n}^{F}(x)-  b_{n}^{F}\right) \frac{z^{n-1}}{F_{n}!} . \label{usefulformulaforsemiclassicalexpansion}
\eea
From this expansion, for $x=1$,
\bea
\frac{e_{F}^{z}-1}{e_{F}^{z}-1}=1+\sum_{n=2}^\infty \left( B_{n}^{F}(1)-  b_{n}^{F}\right) \frac{z^{n-1}}{F_{n}!} \nonumber \\
\Rightarrow \,\,\sum_{n=2}^\infty \left( B_{n}^{F}(1)-  b_{n}^{F}\right) \frac{z^{n-1}}{F_{n}!}=0. \nonumber
\eea
Then, the coefficients at every power of z are zero and
\bea
B_{n}^{F}(1)-  b_{n}^{F}=0,
\eea
for $n=2,3,4,\ldots$
\end{prf}

\section{Acknowledgements}  This work is supported by TUBITAK grant 116F206.

\section{Appendix}
\subsection{Getting $B_{n}^{F}(x)+F_n x^{n-1}$ in two ways}
Our aim is to get equation;
\begin{eqnarray}
B_{n}^{F}(x)+F_n x^{n-1}=H_{n}(x)=\sum_{k=0}^{n} {n \brack k}_{F} B_{n-k}^{F}(x).\nonumber
\end{eqnarray}
Starting with,
\bea
B_{n}^{F}(x)+F_n x^{n-1}&=&H_{n}(x) \Rightarrow \nonumber \\
\sum_{n=0}^{\infty}  B_{n}^{F}(x) \frac{z^n}{F_{n} !}+\sum_{n=0}^{\infty}  F_n x^{n-1} \frac{z^n}{F_{n} !}&=&\sum_{n=0}^{\infty}  H_{n}(x) \frac{z^n}{F_{n} !}\Rightarrow \nonumber \\
\sum_{n=0}^{\infty}  H_{n}(x) \frac{z^n}{F_{n} !} - \sum_{n=0}^{\infty}  B_{n}^{F}(x) \frac{z^n}{F_{n} !}&=&\sum_{n=0}^{\infty}  F_n x^{n-1} \frac{z^n}{F_{n} !}. \nonumber
\eea
For the r.h.s., we get
\bea
z \phantom{.} e^{z x}_{F}=D^{x}_{F} (e^{z x}_{F})=D^{x}_{F} \left(  \sum_{n=1}^{\infty} \frac{x^n z^n}{F_{n} !} \right)= \sum_{n=1}^{\infty}   \frac{F_n \phantom{.} x^{n-1} z^n}{F_{n} !}=\sum_{n=0}^{\infty}   \frac{F_n x^{n-1} z^n}{F_{n} !} .   \nonumber
\eea
Then, we have;
\bea
\sum_{n=0}^{\infty}  H_{n}(x) \frac{z^n}{F_{n} !} - \sum_{n=0}^{\infty}  B_{n}^{F}(x) \frac{z^n}{F_{n} !}&=&z \phantom{.} e^{z x}_{F} \nonumber \Rightarrow \\
\sum_{n=0}^{\infty}  H_{n}(x) \frac{z^n}{F_{n} !} &=&z \phantom{.} e^{z x}_{F}+\sum_{n=0}^{\infty}  B_{n}^{F}(x) \frac{z^n}{F_{n} !} \Rightarrow\nonumber
\eea
\bea
\sum_{n=0}^{\infty}  H_{n}(x) \frac{z^n}{F_{n} !} &=&z \phantom{.} e^{z x}_{F}+\frac{z \phantom{.} e^{z x}_{F} }{e^{z}_{F}-1}\Rightarrow \nonumber \\
\sum_{n=0}^{\infty}  H_{n}(x) \frac{z^n}{F_{n} !} &=&z \phantom{.} e^{z x}_{F}\left(1+\frac{1}{e^{z}_{F}-1}  \right)= z \phantom{.} e^{z x}_{F} \frac{e^{z}_{F}}{e^{z}_{F}-1}\Rightarrow\nonumber \\
\sum_{n=0}^{\infty}  H_{n}(x) \frac{z^n}{F_{n} !} &=&\sum_{n=0}^{\infty}  B_{n}^{F}(x) \frac{z^n}{F_{n} !} \phantom{.} \cdot e^{z}_{F} \Rightarrow\nonumber \\
\sum_{n=0}^{\infty}  H_{n}(x) \frac{z^n}{F_{n} !} &=& \sum_{n=0}^{\infty}  B_{n}^{F}(x) \frac{z^n}{F_{n} !} \phantom{.} \cdot \sum_{k=0}^{\infty} \frac{z^k}{F_{k} !}\Rightarrow \nonumber \\
\sum_{n=0}^{\infty}  H_{n}(x) \frac{z^n}{F_{n} !} &=& \sum_{n=0}^{\infty} \sum_{k=0}^{\infty} \frac{B_{n}^{F}(x)}{F_{n} !} \phantom{.}\frac{z^{n+k}}{F_{k} !} \Rightarrow \nonumber \\
\sum_{n=0}^{\infty}  H_{n}(x) \frac{z^n}{F_{n} !} &{(n+k=N)}{=}& \sum_{N=0}^{\infty} \sum_{k=0}^{N} \frac{B_{N-k}^{F}(x)}{F_{N-k}!} \phantom{.}\frac{z^N}{F_{k} !} \left(  \frac{F_{N}!}{F_{N}!} \right)\Rightarrow \nonumber \\
\sum_{n=0}^{\infty}  H_{n}(x) \frac{z^n}{F_{n} !} &=& \sum_{N=0}^{\infty} \frac{z_{N}}{F_{N}!} \sum_{k=0}^{N} \frac{F_{N}!}{F_{N-k}! F_{k}!} B_{N-k}^{F}(x) \Rightarrow\nonumber \\
\sum_{n=0}^{\infty}  H_{n}(x) \frac{z^n}{F_{n} !} &=& \sum_{N=0}^{\infty} \frac{z_{N}}{F_{N}!} \sum_{k=0}^{N}  {N \brack k}_{F} B_{N-k}^{F}(x) \Rightarrow\nonumber \\
\sum_{n=0}^{\infty}  H_{n}(x) \frac{z^n}{F_{n} !} &{(N=n)}{=}& \sum_{n=0}^{\infty} \frac{z_{n}}{F_{n}!} \sum_{k=0}^{n}  {n \brack k}_{F} B_{n-k}^{F}(x). \nonumber
\eea
By equating two series we get,
\bea
H_{n}(x)=\sum_{k=0}^{n}  {n \brack k}_{F} B_{n-k}^{F}(x).
\eea
Following in another way we show that,
\bea
B_{n}^{F}(x)+F_n x^{n-1}=H_{n}(x)=x^n+ \sum_{j=2}^{n} {n \brack j}_{F} b_{j}^{F} x^{n-j}. \nonumber
\eea
Starting with,
\bea
B_{n}^{F}(x)+F_n x^{n-1}&{(\ref{alternativerepresentationofb.f.p.})}{=}&\sum_{j=0}^{n} {n \brack j}_{F} b_{j}^{F} x^{n-j}+F_n x^{n-1} \nonumber \\
&=&{n \brack 0}_{F}  b_{0}^{F} x^{n}+{n \brack 1}_{F}  b_{1}^{F} x^{n-1}+\sum_{j=2}^{n} {n \brack j}_{F} b_{j}^{F} x^{n-j}+F_n x^{n-1} \nonumber
\eea
\bea
&=& x^n-\frac{F_{n}!}{F_{n-1}! F_{1}!} x^{n-1} +\sum_{j=2}^{n} {n \brack j}_{F} b_{j}^{F} x^{n-j}+F_n x^{n-1} \nonumber \\
&=& x^n +\sum_{j=2}^{n} {n \brack j}_{F} b_{j}^{F} x^{n-j} ,\nonumber
\eea
thus, we obtain
\bea
B_{n}^{F}(x)+F_n x^{n-1}=H_{n}(x)=x^n+ \sum_{j=2}^{n} {n \brack j}_{F} b_{j}^{F} x^{n-j}.
\eea

\end{document}